\documentclass[12pt,oneside,reqno]{amsart}
\usepackage{caption,subcaption}
\usepackage{dcolumn}
\usepackage{amsmath,amsthm,amssymb}
\usepackage{todonotes}
\usepackage{url}
\usepackage{placeins}

\usepackage{hyperref}

\usepackage{cases}

\usepackage{array,booktabs}

\def\be{\begin{equation}}
\def\ee{\end{equation}}

\newtheorem*{theorem*}{Theorem}
\newtheorem*{conj*}{Conjecture}

\theoremstyle{remark}
\newtheorem*{rem*}{Remark}

\def\ra{\rightarrow}

\usepackage{pdflscape}

\usepackage[color, final]{showkeys}
\usepackage{xcolor}

\usepackage{lmodern}
\makeatletter
\ifcase \@ptsize \relax
  \newcommand{\miniscule}{\@setfontsize\miniscule{4}{5}}
\or
  \newcommand{\miniscule}{\@setfontsize\miniscule{5}{6}}
\or
  \newcommand{\miniscule}{\@setfontsize\miniscule{5}{6}}
\fi
\makeatother

\newcounter{todocounter}

\begin{document}

\title[Wave maps on a wormhole]{Multi\-soliton solutions for equivariant\\wave maps on a $2+1$ dimensional wormhole}

\author{Piotr Bizo\'n}
\address{Institute of Theoretical Physics, Jagiellonian
University, 30-348 Krak\'ow, Poland}
\email{piotr.bizon@uj.edu.pl}
\author{Jacek Jendrej}
\address{Sorbonne Universit\'e, Institut de Math\'ematiques de Jussieu -- Paris Rive Gauche, 4 place Jussieu, 75005 Paris, France}
\email{jendrej@imj-prg.fr}
\author{Maciej Maliborski}
\address{University of Vienna, Faculty of Mathematics, Oskar-Morgenstern-Platz~1, 1090 Vienna, Austria, and University of Vienna, Gravitational Physics, Boltzmanngasse 5, 1090 Vienna, Austria}
  \email{maciej.maliborski@univie.ac.at}

\date{\today}
\begin{abstract}
We study equivariant wave maps from the $2+1$ dimensional wormhole to the 2-sphere. This model has  explicit harmonic map  solutions which, in suitable coordinates,  have the form of the sine-Gordon kinks/anti\-kinks.  We conjecture that there exist  asymptotically static chains of $N\geq 2$ alternating  kinks and  anti\-kinks whose subsequent rates of expansion increase in geometric progression as $t\rightarrow \infty$. Our argument employs the method of collective coordinates  to derive effective finite-dimensional ODE models for the asymptotic dynamics of $N$-chains. For $N=2,3$ the predictions of these effective models are verified  by direct PDE computations which demonstrate that the $N$-chains  lie at the threshold of kink-anti\-kink annihilation.

\end{abstract}
\maketitle

\section{Introduction}
For nonlinear evolution equations the mapping from the set of initial data at $t=0$ to the set of end states at $t=\infty$ (or $t=T<\infty$ in the case of blowup at time $T$) is typically highly surjective meaning that the set of end states is, in a sense, small. In particular, for nonlinear wave equations on unbounded domains  this property of dynamical asymptotic simplification  is at the heart of the soliton resolution conjecture which states that  solutions decompose asymptotically into a finite superposition of
decoupled  bound states (solitons) and free outgoing radiation \cite{T}.

  In recent years there has been great progress in mathematical understanding of the soliton resolution for radial solutions of the energy critical nonlinear wave equations (see \cite{K} for a survey and references therein). Notably, the proof of the soliton resolution involving complete taxonomy of possible end states of evolution was given for
equivariant wave maps $\mathbb{R}^{2+1}\ra \mathbb{S}^2$ \cite{DKMM, JL} and the equivariant Yang-Mills (YM) equation in $4+1$ dimensions \cite{JL}.
We remark that these results are abstract in the sense that they do not say which of possible end states  are actually realized in the evolution.

While most of the mathematical work on soliton resolution conjecture has been concerned with nonlinear radial wave equations on $d+1$ Minkowski spacetime, it was pointed out by the first author \cite{B1} that the problem becomes more tractable  for the corresponding equations posed  on a wormhole spacetime. 
Subsequently, the soliton resolution conjecture for equivariant wave maps from the $3+1$ dimensional wormhole  into the 3-sphere was formulated and verified numerically in \cite{BK} and later  proved by Rodriguez \cite{R1,R2}. The simplicity of this model stems from the fact that in each topological sector there is a unique endstate of evolution given by the rigid stable soliton.
In this paper we consider the analogous model in two spatial dimensions, i.e.
equivariant wave maps from the $2+~1$ dimensional wormhole  into the 2-sphere. This model admits  an explicit  static solution (harmonic map of degree one, hereafter called the kink), whose reflections (called anti\-kinks) and conformal transformations serve as building blocks for a variety of multi\-soliton solutions. Here, as the first step toward the taxonomy of possible end states, we consider multi\-soliton solutions of degree zero and one which are asymptotically static. Although such ``parabolic" solutions are not generic, they play an important role as separatrices between different generic behaviors in the global phase portrait. We note that similar solutions in the case of the $\phi^{4}$ model were considered in \cite{JKL1, JL9}.

The paper is organized as follows. In Sec.~2 we introduce the model and formulate our main conjecture asserting the existence of asymptotically static chains of $N$ alternating kinks and anti\-kinks whose positions and  mutual distances tend to infinity  as $t\rightarrow \infty$. In Sec.~3 we derive an effective finite dimensional system of ordinary differential equations that approximate the asymptotic dynamics of $N$-chains. Finally, in Sec.~4 we present the results of direct numerical computations of the wave map equation and confirm the predictions of the ODE model for $N=2,3$.

\section{Model}

A map $X: M\mapsto N$ from a Lorentzian manifold $(M,g_{\alpha\beta})$ into a Riemannian manifold $(N,G_{AB})$
is a wave map if it is a critical point of the action
\be
S[X]=\int_M g^{\alpha\beta} \partial_{\alpha} X^A \partial_{\beta} X^B G_{AB},
\ee
which gives the wave map equation
\be \label{wmeq}
\Box_g X^A + \Gamma_{BC}^A(X) \partial_{\alpha}X^B \partial_{\beta}X^C g^{\alpha\beta}=0,
\ee
where $\Gamma_{BC}^A$ are the Christoffel symbols of $G_{AB}$.
\vskip 0.1cm
In this paper, as domain we take the three-dimensional Lorentzian manifold $M=\{(t,r)\in \mathbb{R}^2, \phi\in \mathbb{S}^1\}$ with the metric
\be \label{worm_metric}
    ds^2=-dt^2 + dr^2 + (r^2+a^2) d\phi^2,
\ee
where $a$ is a positive constant. This ultrastatic spacetime is the prototype example of a wormhole geometry which has two asymptotically flat ends at $r \rightarrow \pm \infty$ connected by the circular neck   of circumference $2\pi a$ at $r=0$
(see \cite{MT} for a pedagogical introduction to the geometry and physics of wormholes).

As the target we take the unit 2-sphere $N=\mathbb{S}^2$ with the round metric $ds^2=dU^2+\sin^2{U}\, d\Phi^2$.
We assume the equivariant ansatz:
  $U=U(t,r)$ and $\Phi=k \phi$, where $k$ is a positive integer. Then, the wave map equation \eqref{wmeq} reduces to the semilinear scalar wave equation
  \be \label{eqr}
  U_{tt} = U_{rr}+\frac{r}{r^2+a^2} \, U_r -\frac{k^2}{2} \,\frac{\sin(2U)}{r^2+a^2}\,.
  \ee
The length scale $a$ removes the singularity at $r=0$, thereby ensuring global-in-time regularity of solutions.  In the following we set $a=1$ by the choice of  unit of length and  introduce a new coordinate $x=\sinh^{-1}(r)$. In terms of~$x$ the spatial part of the wormhole metric  is manifestly conformal to the flat metric on the cylinder $\mathbb{R}\times \mathbb{S}^1$,
\begin{equation}\label{worm_x}
ds^2=-dt^2+\cosh^2{\!x} \,(dx^2+d\phi^2),
\end{equation}
and Eq.~\eqref{eqr} for $u(t,x)=U(t,r)$ simplifies  to
\be \label{eqx}
\cosh^2{\!x}\, u_{tt} = u_{xx} -\frac{k^2}{2}\,\sin(2u).
\ee
The conserved energy is $E=T+V$, where
\be
 \label{energy}
T = \frac{1}{2}\, \int\limits_{-\infty}^{\infty} \cosh^2{\!x}\, u_t^2 \, dx, \qquad
V = \frac{1}{2}\, \int\limits_{-\infty}^{\infty} \left( u_x^2  +k^2\,\sin^2{\!u}\right) \, dx
\ee
are the kinetic and potential energies.
Any finite energy solution must be equal to integer multiples of $\pi$ at $x=\pm \infty$.
Without loss of generality we set $u(t,-\infty)=0$, so that the integer $n=u(t,+\infty)/\pi$ determines the topological degree of solution (which is preserved in evolution).
\vskip 0.1cm
For $n=0$ the energy attains the minimum $E=0$ at the vacuum $u=0$.  For  $n=1$  we have the Bogomol'nyi inequality
\be \label{bogom}
E(u)\geq V(u)=\frac{1}{2}\int_{-\infty}^{\infty}
\left(u_x-k \sin{u}\right)^2\, dx +2k \geq 2k,
\ee
which is saturated on the kink solution
\be
u=Q(x)=2 \arctan\left(e^{kx}\right).
\ee
For static solutions equation \eqref{eqx} is  invariant under  translations in $x$ (which is due to the conformal invariance of  harmonic maps in two dimensions), hence translations of the kink $Q(x)$ generate a  one-parameter family  of kinks $Q_c(x):=Q(x-c)$ with the same energy $2k$. By reflection symmetry, there is the corresponding family of
 anti\-kinks $-Q_c(x)$ of degree $n=-1$. For $|n|\geq 2$ there is no static solution and $E>2k|n|$, where the lower bound is attained by an ideal configuration of $n$ kinks (or anti\-kinks) separated by infinite mutual distances.
 \vskip 0.1cm
Next, we briefly discuss stability of the kink under small perturbations.
  The linearized operator  about the kink
  \begin{equation}\label{L}
   L=-\frac{1}{\cosh^2{x}}\, \left(\partial_{xx}-k^2 \cos(2Q_c(x)\right)
 \end{equation}
 is self-adjoint with respect to the scalar product
$(f,g)=\int_{-\infty}^{\infty} f(x) g(x) \cosh^2{\!x}\, dx$.  The continuous spectrum is $[0,\infty)$, where $0$ at the bottom of the spectrum is associated with the zero mode given by the vector tangent to the curve $Q_c(x)$
\begin{equation}\label{0mode}
\xi(x):= -\frac{dQ_c(x)}{dc} =\frac{k}{\cosh{k(x-c)}}.
\end{equation}
For $k=1$ the  zero mode is not square-integrable, while  for $k\geq 2$ it is a genuine eigenfunction.  Since the zero mode has no zeroes, it follows from the Sturm oscillation theorem that there are no negative eigenvalues.
To avoid technical difficulties, hereafter we assume that $k\geq 2$.
\vskip 0.1cm
  Using the fact that $Q_{c}$ is a degenerate minimizer of energy,
  it is routine to show that $Q_{c}$ is orbitally stable, i.e. if at $t=0$ the solution is close to $Q_{c(0)}$ for some $c(0)$, then it remains close to the orbit $\mathcal{O}_Q:=\{Q_{c}: c\in \mathbb{R}\}$ for all~$t$.
  Moreover,  the kink is  expected to be  asymptotically orbitally stable, i.e. if  the solution is initially close to $\mathcal{O}_Q$, then it converges (locally in space) to $\mathcal{O}_Q$ as $t\rightarrow \infty$. To establish this, one needs to analyze  interaction of the moving kink with radiation. We leave  this difficult open problem for future investigations and focus here on special  multi\-kink solutions for which radiation is asymptotically negligible. Namely,  we shall give analytic and numerical evidence for the existence
  of asymptotically static $N$-chains of alternating kinks and anti\-kinks which, up to an error which tends to zero, are given  by
\begin{subnumcases}
  {\label{Nchain} u_N(t,x):=}
  \sum\limits_{j=1}^J (-1)^j \left[Q(x+c_j(t))- Q(x-c_j(t))\right], & \!\!\!\!\!\!\!\!\!\!\!\!\!\!\!\!\!\!\!\!\!\!\!\!\!\!\!\!\!\! $N=2J$,\\
  Q(x)+\sum\limits_{j=1}^J (-1)^j \left[Q(x+c_j(t))+ Q(x-c_j(t))\right], &  \nonumber\\
  & \!\!\!\!\!\!\!\!\!\!\!\!\!\!\!\!\!\!\!\!\!\!\!\!\!\!\!\!\!\! $N=2J+1$
\end{subnumcases}
where $c_j(t)\rightarrow \infty$ and $c_{j}(t)-c_{j-1}(t) \rightarrow \infty$ for all $1\leq j\leq J$ as $t\rightarrow \infty$, where $c_0=0$. Note that $\lim_{t\rightarrow \infty}E(u_N)=2kN$.
\vskip 0.1cm
   For simplicity of the presentation, in what follows we set the equivariance index $k$ equal to 2. The results for $k>2$ are qualitatively the same, while the case $k=1$ is left as a challenge.
\section{ODE approximation}
In this section we derive   finite-dimensional ODE models for the asymptotic dynamics of the $N$-chains \eqref{Nchain}  via the method of collective coordinates. This method works as follows;  see \cite{MS} for the general exposition.  The Lagrangian functional $\mathcal{L}=T-V$ for the wave map equation \eqref{eqx} is evaluated on the $N$-chain configuration $u_N(t,x)$ by integrating over $x$. The resulting function  is interpreted as an effective  Lagrangian $\mathcal{L}_{\text{eff}}$ which determines the evolution of the kink and anti\-kink positions $c_j(t)$. Thus,
 instead of extremizing $\mathcal{L}$ with respect to all variations of $u(t,x)$, one considers only variations allowed by the ansatz.
This approach omits radiative degrees of freedom so it is expected to provide good approximation to true dynamics only if radiation is negligible, which is precisely the case at hand.

Let us first consider the 2-chain
$
u_2(t,x)=Q(x+c(t))-Q(x-c(t))
$
for late times and hence for large $c(t)$.
Evaluating the kinetic and potential energies for $k=2$, we get for large $c$
\be \label{Lu2}
T(u_2)\approx \frac{\pi}{2}  e^{2 c} {\dot c}^2,\qquad
V(u_2)\approx 8-16 e^{-4c}\,,
\ee
from which we see that the effective kinetic energy of an isolated kink positioned at some large $c(t)$ is approximately equal to $\frac{\pi}{4} e^{2 c} {\dot c}^2$, while the binding energy of the  kink-anti\-kink pair separated by a large distance $d$ is approximately equal to $-16 e^{-2d}$. Using this, we get the effective Lagrangians  which determine the asymptotic dynamics of $N$-chains (here by convention $c_0=0$)
\begin{subnumcases} {\label{Leff} \mathcal{L}_{\text{eff}}=}
\frac{\pi}{2} \sum\limits_{j=1}^J e^{2 c_j} {\dot c}_j^2 + 32 \sum\limits_{j=1}^J e^{-2(c_j-c_{j-1})}, & \!\!\!\!\!$N=2J+1$,\\
\frac{\pi}{2}  \sum\limits_{j=1}^J e^{2 c_j} {\dot c}_j^2 + 16 e^{-4 c_1} + 32 \sum\limits_{j=2}^J e^{-2(c_j-c_{j-1})}, & \!\!\!\!\!$N=2J$.
\end{subnumcases}
It is convenient to define $r_j=e^{c_j}$ and rescale time $t\rightarrow \frac{8}{\sqrt{\pi}} t$. Then the effective Lagrangians read (where  $r_0=1$)
\begin{subnumcases} {\label{Leff} \mathcal{L}_{\text{eff}}=}
\sum\limits_{j=1}^J  {\dot r}_j^2 +  \sum\limits_{j=1}^J \frac{r_{j-1}^2}{r_j^2}, & $N=2J+1$, \label{lag_odd}\\
\sum\limits_{j=1}^J {\dot r}_j^2 + \frac{1}{2 r_1^4} + \sum\limits_{j=2}^J \frac{r_{j-1}^2}{r_j^2}, & $N=2J$,\label{lag_even}
\end{subnumcases}
which can be interpreted as Lagrangians for $J$ point particles on half-line\footnote{One can arrive at \eqref{Leff} directly from the wave map Lagrangian in the original coordinates $(t,r)$.  In this formulation it is natural to view $r_j(t)$ as expanding concentric circles.} with rather unusual attractive nearest-neighbor interactions. We emphasize that these effective Lagrangians are valid only in the asymptotic regime $r_j\rightarrow \infty$ and $r_{j}/r_{j-1} \rightarrow \infty$. Next, we analyze the odd and even $N$-chains  separately.
\begin{center}
$N=2J+1$
\end{center}
The Euler-Lagrange equations corresponding to the Lagrangian \eqref{lag_odd} are (here by convention $r_{J+1}=\infty$)
\begin{equation}\label{eq_odd}
  \ddot r_{j} = -\frac{r^2_{j-1}}{r^3_j} + \frac{r_{j}}{r^2_{j+1}},\qquad j=1,\dots, J.
\end{equation}
and the associated  conserved energy is
\begin{equation}\label{E_odd}
  E_{\text{eff}} =\sum_{j=1}^J \dot{r_j}^2 -\sum_{j=1}^J \frac{r_{j-1}^2}{r^2_{j}},
\end{equation}
where by convention the effective potential energy is measured with respect to the limit potential energy $\lim_{t\rightarrow \infty} V(u_N)=4N$, hence the effective evolution is constrained to the zero energy level $E_{\text{eff}}=0$. We claim that there exists a zero-energy asymptotically static solution of the system \eqref{eq_odd}
such that for $t\ra\infty$
\begin{equation}\label{ss_odd}
  r_j(t) \sim (A t)^{\frac{j}{J+1}}, \quad A=\frac{J+1}{\sqrt{J}}.
\end{equation}
We shall denote this solution by $r_j^{(N)}(t)$.
For  $J=1$ we have
\be \label{r3}
 r_1^{(3)}(t)=(2\, t)^{\frac{1}{2}},
  \ee
  which is exact. For $J\geq 2$, the solution $r_j^{(2J+1)}(t)$ can be constructed formally as follows. Assuming temporarily that the right-hand sides of equations \eqref{eq_odd} vanish for $1\leq j\leq J-1$, one gets $r_{j-1} r_{j+1}=r_j^2$, thus  $r_j$ form a geometric sequence
\begin{equation}\label{rj}
   r_{j}=r_1^j,
\end{equation}
and therefore $r_{J-1}=r_J^{\frac{J-1}{J}}$. Substituting this into Eq.~\eqref{eq_odd} for $j=J$ gives
$
  \ddot r_J = - r_J^{-\frac{J+2}{J}},
$
whose zero energy solution is
$
  r_J(t)= (A t)^{\frac{J}{J+1}}$.
Combining this  result with \eqref{rj} yields  \eqref{ss_odd}. Having the leading order asymptotics, the solution can be systematically improved by adding subleading corrections. Let us illustrate how this  goes for $J=2$, where \eqref{ss_odd} reads
\begin{equation}\label{r2_odd}
  r_j(t) \sim (At)^{\frac{j}{3}},  \quad A=\frac{3}{\sqrt{2}}.
\end{equation}
To derive the asymptotic series expansion for $t\ra \infty$ starting from \eqref{r2_odd}, we introduce a new independent variable $\tau=(At)^{1/3}$ and  factor out the leading asymptotics $r_j(t)=\tau^j v_j(\tau)$.
Substituting this into the system \eqref{eq_odd} we get
\begin{subequations}\label{eqv5}
\begin{eqnarray}
&&\tau^2 \ddot{v}_1 - 2 v_1 +2\tau^2 \left(\frac{1}{v_1^3}-\frac{v_1}{v_2^2}\right)=0,\\
&&\tau^2 \ddot{v}_2  +2 \tau \dot{v}_2 -2 v_2 +\frac{2v_1^2}{v_2^3}=0.
\end{eqnarray}
\end{subequations}
Starting from $v_1=v_2=1$ and repeatedly using  the method of dominant balance we find the formal solution given by the following asymptotic series for $\tau\rightarrow\infty$
\begin{subequations}\label{v5_as}
\begin{eqnarray}
\hspace{-1.5cm}&& v_1(\tau)=1-\frac{3}{8}\,\tau^{-2}+c\tau^{-3}+\frac{3}{128}\,\tau^{-4}+\frac{3c}{8}\,\tau^{-5}+ \left(\frac{443}{3072}-c^2\right)\,\tau^{-6}+...,\\
\hspace{-1.5cm}&& v_2(\tau)=1-\frac{1}{4}\tau^{-2}+2c\tau^{-3}-\left(\frac{5}{192}+c^2\right)\,\tau^{-6}+...
\end{eqnarray}
\end{subequations}
where the free parameter $c$ is the trace of time translation symmetry of the original system. This concludes the derivation of the zero-energy asymptotically static solution of the system \eqref{eq_odd} for $J=2$
\be \label{chain5}
r^{(5)}_j(t)=(A t)^{\frac{j}{3}}\, v_j\left((At)^{\frac{1}{3}}\right),
\ee
 where $v_j(\tau)$ are given in \eqref{v5_as}.
\vskip 0.1cm
Next, we show that the solution $r^{(5)}_j(t)$ has one unstable mode \emph{within} the zero-energy level.
To this end, we  insert\footnote{For the purpose of this calculation the subleading terms in the asymptotic expansions \eqref{v5_as} can be omitted.} $v_j(\tau)=1+\varepsilon \xi_j(\tau)$ into the system~\eqref{eqv5}. Keeping only linear terms in $\varepsilon$, we get
\begin{subequations}\label{eq5_xi}
\begin{eqnarray}
&&\tau^2 \ddot{\xi_1}  - 2\xi_1 +4\tau^2 (\xi_2-2 \xi_1)=0,\\
&& \tau^2 \ddot{\xi_2} + 2\tau \dot{\xi}_2 -8 \xi_2 + 4 \xi_1 =0.
\end{eqnarray}
\end{subequations}
Using the method of dominant balance\footnote{Alternatively,  one can solve the system \eqref{eq5_xi}  exactly in terms of generalized hypergeometric functions and then take  their asymptotic expansions at infinity}  we obtain the asymptotic behavior of growing mode for $\tau \ra\infty$
\begin{equation}\label{xi5}
\xi_1(\tau)\sim  e^{\sqrt{8}\, \tau}, \qquad \xi_2(\tau)\sim -\frac{1}{2} \tau^{-2}   e^{\sqrt{8}\, \tau}\,.
\end{equation}
For the  perturbation generated by this mode $\delta E_{\text{eff}}=\mathcal{O}(\varepsilon^2)$,   hence it is tangent to  the zero-energy surface, as claimed. Notice the opposite signs of $\xi_1(\tau)$ and $\xi_2(\tau)$. In addition, there is a polynomially growing unstable mode $\xi_1 \sim \varepsilon \tau^2, \xi_2 \sim 2\varepsilon \tau^2$ for which $\delta E_{\text{eff}} =\mathcal{O}(\varepsilon)$.
\newpage
\begin{center}
$N=2J$
\end{center}
The Euler-Lagrange equations corresponding to the Lagrangian \eqref{lag_even} are
\begin{equation}\label{eq_even}
  \ddot r_1 = -\frac{1}{r_1^5}+\frac{r_1}{r_2^2} \quad \mbox{and}\quad \ddot r_{j} = -\frac{r^2_{j-1}}{r^3_j} + \frac{r_{j}}{r^2_{j+1}} \quad\mbox{for}\quad  j=2,\dots, J.
\end{equation}
By an analogous calculation as in the case of odd $N$, we obtain a zero-energy asymptotically static solution of the system \eqref{eq_even}
such that for $t\ra\infty$
\begin{equation}\label{r_even}
  r_j(t) \sim (A t)^{\frac{2j-1}{2J+1}}, \quad A=\frac{2J+1}{\sqrt{4J-2}}.
\end{equation}
For $J=1$ the solution
\be \label{r2}
 r_1^{(2)}(t)=\left(\frac{3}{\sqrt{2}}\, t\right)^{\frac{1}{3}},
  \ee
 is exact, while for $J\geq 2$  the leading order asymptotics can be systematically improved. As before, we give the details in the case $J=2$, where
\begin{equation}\label{r2_even}
 r_j(t)\sim (At)^{\frac{2j-1}{5}},\qquad  A=\frac{5}{\sqrt{6}}.
\end{equation}
We define a new independent variable $\tau=(At)^{1/5}$ and  factor out the leading asymptotics $r_j(t)=\tau^{2j-1} v_j(\tau)$.
Substituting this into the system \eqref{eq_even} we get
\begin{subequations}\label{eqv4}
\begin{eqnarray}
 && \tau^2 \ddot{v}_1 - 2\tau \dot{v}_1 -4 v_1 +6 \tau^4 \left(\frac{1}{v_1^5}-\frac{v_1}{v_2^2}\right)=0,\\
 && \tau^2 \ddot{v}_2 +2 \tau \dot{v}_2 +6 \left(\frac{v_1^2}{v_2^3}-v_2\right)=0.
  \end{eqnarray}
\end{subequations}
Starting from $v_1=v_2=1$ and repeatedly using  the method of dominant balance we find the formal solution given by the following asymptotic series for $\tau\rightarrow\infty$
\begin{subequations}\label{v4}
\begin{eqnarray}
&& v_1 = 1-\frac{1}{6} \tau^{-4} + \frac{c}{3} \tau^{-5} -\frac{83}{1944} \tau^{-8}...,\\
&&  v_2 = 1-\frac{1}{6}\tau^{-4} + c\tau^{-5} +\frac{7}{648} \tau^{-8}+...
\end{eqnarray}
\end{subequations}
Thus,  the zero-energy asymptotically static solution of the system \eqref{eq_even} for $J=2$  is
\be
r^{(4)}_j(t)=(A t)^{\frac{2j-1}{5}}\, v_j\left((At)^{\frac{1}{5}}\right),
\ee
 where $v_j(\tau)$ are given in \eqref{v4}.

To find the growing mode about this solution we insert $v_j(\tau)=1+\varepsilon \xi_j(\tau)$ into the system \eqref{eqv4}. Keeping only linear terms in $\varepsilon$, we get
\begin{subequations}
\begin{eqnarray}
&&\tau^2 \ddot{\xi_1}  -2 \tau \dot{\xi}_1 - 4\xi_1 +12 \tau^4 \left(\xi_2-3\xi_1\right)=0,\\
&& \tau^2 \ddot{\xi_2} + 2\tau \dot{\xi}_2 - 24 \xi_2 + 12\xi_1 =0.
\end{eqnarray}
\end{subequations}
From this we obtain the growing mode for $\tau\ra\infty$
\begin{equation} \label{xi4}
 \xi_1(\tau)\sim  \tau^{\frac{1}{2}} e^{3\tau^2}, \qquad
 \xi_2(\tau)\sim -\frac{1}{3} \tau^{-\frac{7}{2}}  e^{3\tau^2}\,.
\end{equation}
As before, one can check  that this mode is tangent to the zero energy surface. In addition, there is  a polynomially growing unstable mode $\xi_1 \sim \varepsilon \tau^4, \xi_2 \sim 3\varepsilon \tau^4$ for which $\delta E_{\text{eff}} =\mathcal{O}(\varepsilon)$.
\vskip 0.1cm
\noindent \emph{Remark.} On the basis of odd and even cases with $J=2$ discussed above and similar calculations for $J=3$, we conjecture that the stable manifolds of $2J$ and $(2J+1)$-chains have codimension $J$. For this reason it seems very difficult to verify the above ODE predictions for $N>3$ (i.e. $J>1$)  by direct PDE computations\footnote{Our unsuccessful attempts to achieve this goal delayed this paper significantly.}  and consequently in the next section we present results only for $N=2,3$ (i.e. $J=1$).
\section{PDE computations}
In this section we solve Eq.~\eqref{eqx} numerically using the hyperboloidal formulation of the initial value problem \cite{An}. To this end we  define a new time  coordinate $s=t-\cosh{x}$.
 The hypersurfaces of constant~$s$ are ``hyperboloidal'' which means that they are spacelike hypersurfaces that approach the future null infinities of the wormhole spacetime  along the outgoing null cones for $x\rightarrow \pm \infty$.
In terms of $s$ Eq.~\eqref{eqx} becomes (for $k=2$)
\begin{equation}\label{eqxs}
  u_{ss}+2 \sinh(x)\, u_{sx} +\cosh(x)\, u_s= u_{xx}-2 \sin(2u).
\end{equation}
We shall solve this equation for smooth finite energy initial data \linebreak $(u(0,x), u_s(0,x))$ of degree $n=0,1$. From the general theory of wave equations on asymptotically flat spacetimes it is known that  solutions behave as follows at null infinities \cite{CL}
\begin{subnumcases}{\label{exp_scri}
u(s,x)\sim}
b_{-}(s)\, e^{\frac{x}{2}} & $x \rightarrow -\infty$\\
n\pi+b_{+}(s)\, e^{-\frac{x}{2}} & $x \rightarrow \infty$,
\end{subnumcases}
where $b_{\pm}(s)$ are the radiation coefficients.
Multiplying Eq.~\eqref{eqxs} by $u_s$ and rearranging terms we obtain the conservation law
\begin{equation}\label{claw}
\partial_s \left(\frac{1}{2} u_s^2 + \frac{1}{2}  u_x^2 + 2 \sin^2{\!u}\right)
   +  \partial_x \left(\sinh(x)  u_s^2 - u_s \,u_x \right)=0.
\end{equation}
Integrating this over $x$ and using \eqref{exp_scri} we get  the energy loss formula
\begin{equation}\label{loss}
  \frac{d\mathcal{E}}{ds} = -\frac{1}{2}{\dot b}_{-}^2(s)-\frac{1}{2}{\dot b}_{+}^2(s),
\end{equation}
where we defined the Bondi-type energy
\begin{equation}\label{bondi}
  \mathcal{E}(u):=\int_{-\infty}^{\infty} \left(\frac{1}{2}  u_s^2 + \frac{1}{2} u_x^2 + 2\sin^2{\!u}\right) dx.
\end{equation}
Note that the potential part of the energy $\mathcal{E}$ is the same as $V(u)$ for the conserved energy $E$, however the kinetic part is different.
Since the Bondi energy is monotonically decreasing and bounded from below---by $0$ for $n=0$ and by $4$ for $n=1$ [due to Bogomolnyi's inequality \eqref{bogom}]---there exists a limit
\be \label{E_inf}
\mathcal{E}_{\infty}=\lim_{s\rightarrow \infty} \mathcal{E}(u(s)).
\ee
The energy of radiation goes  to zero as $s\rightarrow \infty$,  hence the limiting energy is quantized $\mathcal{E}_{\infty}=4N$, where $N$ is the final number of decoupled kinks and anti\-kinks.

In numerical studies we introduce the coordinate $y=\tanh\left(x/4\right)$ which  compactifies the spatial domain to the interval $[-1,1]$ and transforms Eq.~\eqref{eqxs} to
\begin{equation}\label{eqys}
   u_{ss} +\frac{2y(1+y^2)}{1-y^2}\,u_{sy} +\frac{y^4+6y^2+1}{(1-y^2)^2}  u_s =\frac{1}{16} (1-y^2)\partial_y \left[ (1-y^2) u_y \right]-2 \sin(2u).
\end{equation}
This equation is solved numerically using the method of lines with Chebyshev pseudospectral spatial discretisation.
We work with nodal representation on the Chebyshev-Gauss-Lobatto grid points where the derivatives are replaced by corresponding matrix multiplications of differentiation matrices \cite{Tr}.
Since we consider only odd or even solutions,  computations are done on the half-domain $0\leq y\leq 1$ together with appropriate boundary conditions at $y=0$.
For time integration of the resulting ODEs for function values at grid points we use IDA method \cite{IDA} as implemented in Wolfram \textit{Mathematica} \cite{Mathematica}.
Setting a stringent error tolerance  in the time-stepping algorithm we observe fast spectral convergence with increasing number of grid points.

\begin{figure}[ht]
	\centering
	\includegraphics[width=1.\textwidth]{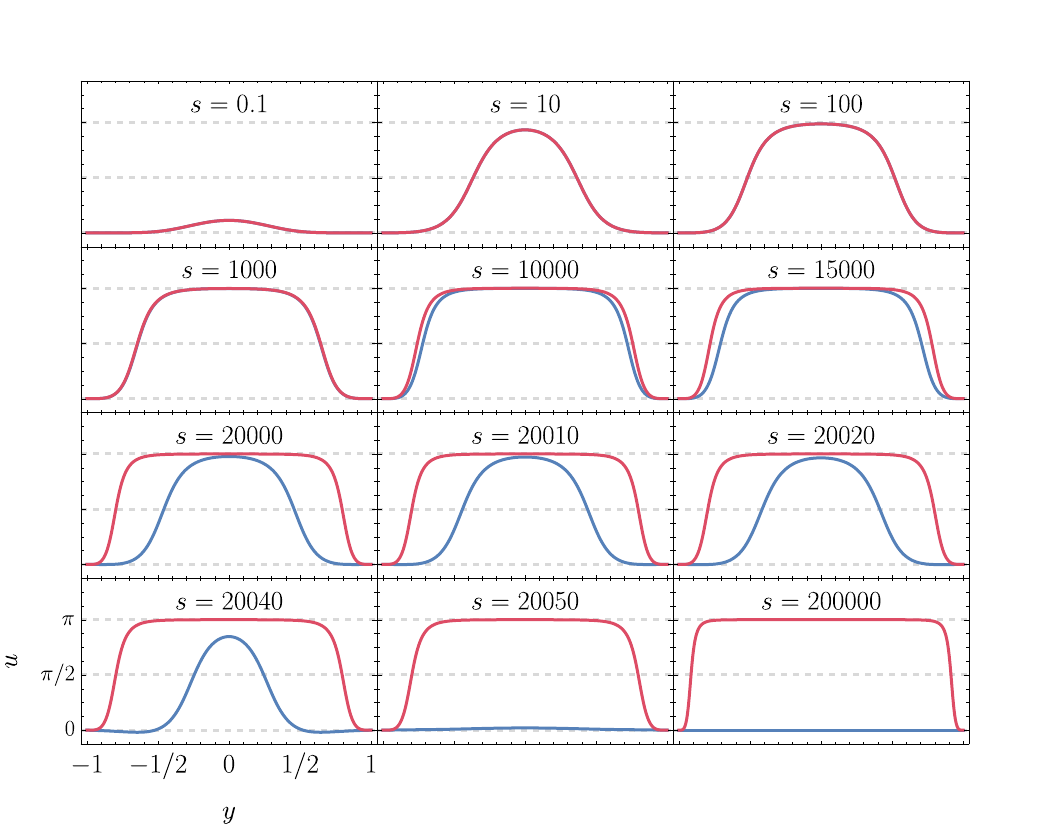}
	\caption{\small{Snapshots of even solutions for marginally  subcritical initial data \eqref{id1} with $b=b_{*} - 6\times10^{-16}$ (in blue) and marginally supercritical data with $b=b_{*} + 6\times10^{-16}$ (in red).\vspace{0cm} }
	}
	\label{fig:N2evolution}
\end{figure}

\begin{figure}[ht]
	\centering
	\includegraphics[width=0.67\textwidth]{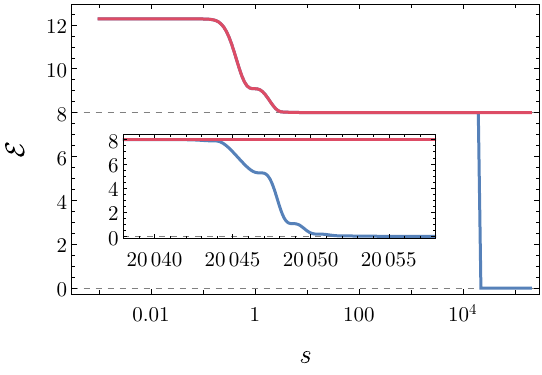}
	\caption{\small{The Bondi energy of sub- and supercritical solutions from Fig.~\ref{fig:N2evolution}. For intermediate times the energy is slightly above the energy $\mathcal{E}=8$ of the asymptotically static kink-anti\-kink pair. For supercritical solutions the energy tends slowly to $\mathcal{E}=8$ for $s\ra \infty$. For subcritical solutions the energy is rapidly radiated away 
 during the collision of kink and anti\-kink and tends to $0$ (see the inset).\vspace{0.5cm} }}
	\label{fig:N2energy}
\end{figure}

\begin{figure}[ht]
	\centering
	\includegraphics[width=0.67\textwidth]{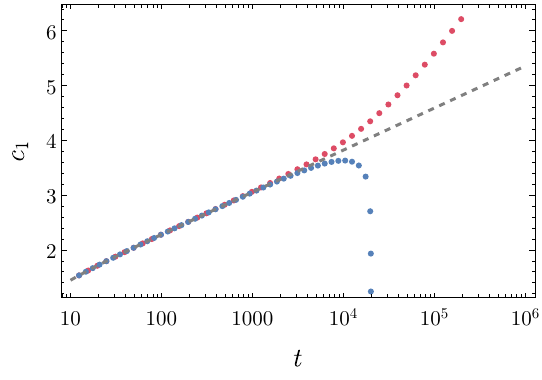}
	\caption{\small{Position of the kink $c_1(t)$ for sub- and supercritical from Fig.~\ref{fig:N2evolution}. The gray dashed curve is the fit of formula \eqref{fitN2} to the numerical data for intermediate times. For later times the solutions deviate from criticality along the unstable directions. }
	}
	\label{fig:N2location}
\end{figure}

We now discuss numerical results for sample even and odd initial conditions. As even initial data we chose
\begin{equation}
	\label{id1}
	u(0,y) = 0\,,
	\quad
	\partial_{s}u(0,y) = b \exp\left(-\frac{4}{(1-y^{2})^{2}}+4\right)
	\,,
\end{equation}
where $b$ is an adjustable parameter.
For small values of $b$ solutions tend to the  trivial map $u=0$.
For larger values of $b$  we observe the formation of kink-anti\-kink pair. The pair is initially expanding and its evolution can be viewed as competition between the kinetic energy and the attractive force between the kink and anti\-kink. For large $b$ the kinetic energy wins and the pair expands to infinity. For smaller $b$ the attraction wins; after some time the expansion stops  and the pair starts contracting which eventually leads to collision and annihilation. 
At the threshold between these two scenarios (for $b_{*}\approx 3.78523$)  we observe the asymptotically static $N=2$ configuration $u_2$ given in (12a) with the rate of expansion as predicted by the ODE approximation \eqref{r2}. To verify that, for intermediate times we fit the position of the anti\-kink $c_1(t)$ to the  formula \eqref{r2} with a free parameter~$A$
\be \label{fitN2}
	c_1(t) = \frac{1}{3}\left(\log(t-t_{0})+\log A\right)
\end{equation}
and get $ A \approx 9.45795$ which is  in  good agreement with the ODE prediction $A=12 \sqrt{\frac{2}{\pi }}\approx 9.57461$ (which follows from \eqref{r2} and  rescaling $t\rightarrow \frac{8}{\sqrt{\pi}} t$).
\vskip 0.1cm
The snapshots of marginally sub- and supercritical solutions are depicted in Fig.~\ref{fig:N2evolution}, while 
the corresponding Bondi energies and positions of the anti\-kink are shown in Fig.~\ref{fig:N2energy} and Fig.~\ref{fig:N2location}


\noindent For odd initial data we chose a one-parameter family
\begin{equation}
	\label{id2}
	u(0,y) = \frac{\pi}{2} \left[\left(1+\sin \left(\frac{\pi  y}{2}\right)\right)-4b
   \sin \left(\frac{\pi  y}{2}\right) \cos ^2\left(\frac{\pi 
   y}{2}\right)\right]
   \,,
   \quad
   \partial_{s}u(0,y) = 0
   \,.
\end{equation}
Here the subcritical solutions converge to the kink, while supercritcal solutions converge to the kink-anti\-kink-kink configuration with the static anti\-kink positioned at the center and outer kinks expanding to infinity.
At the threshold between these two scenarios (for $b_{*}\approx 0.381895$)  we observe the asymptotically static $N=3$ configuration $u_3$ given in (12b) with the rate of expansion as predicted by the ODE approximation \eqref{r3}. To verify that, for intermediate times we fit the position of the kink $c_1(t)$ to the  formula \eqref{r3} with a free parameter~$A$
\be \label{fitN3}
	c_1(t) = \frac{1}{2}\left(\log(t-t_{0})+\log A\right)
\end{equation}
and get $A \approx 9.05616$ which is  in good agreement with the ODE prediction $A=\frac{16}{\sqrt{\pi}}\approx 9.02703$.
\begin{figure}[ht]
	\centering
	\includegraphics[width=0.95\textwidth]{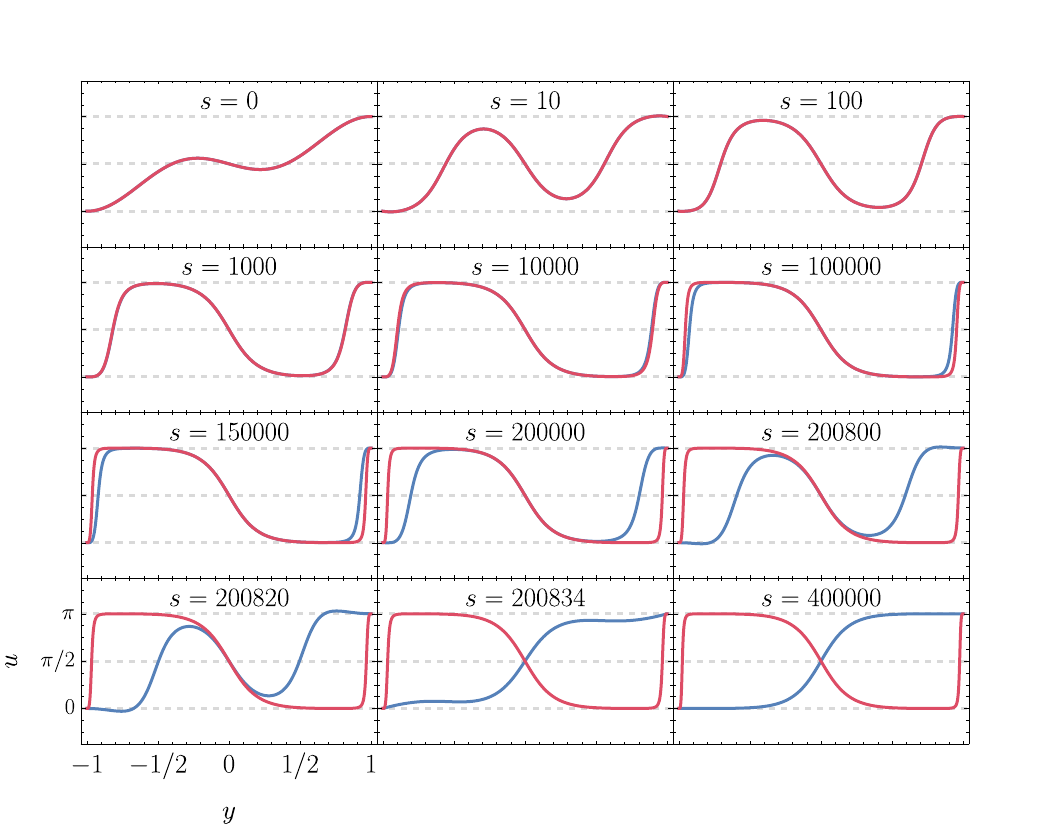}
	\caption{\small{Snapshots of odd solutions for marginally  subcritical initial data \eqref{id2} with $b=b_{*} - 6\times10^{-16}$ (in blue) and marginally supercritical data with $b=b_{*} + 6\times10^{-16}$ (in red). }
	}
	\label{fig:N3evolution}
\end{figure}

\begin{figure}[ht]
	\centering
	\includegraphics[width=0.61\textwidth]{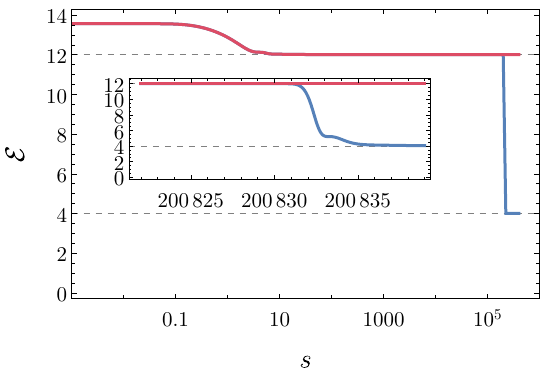}
	\caption{\small{The Bondi energy of sub- and supercritical solutions from Fig.~\ref{fig:N3evolution}. For intermediate times the energy is slightly above the energy $\mathcal{E}=12$ of the asymptotically static kink-anti\-kink-kink configuration. For supercritical solutions the energy tends to $\mathcal{E}=12$ for $s\ra \infty$. For subcritical solutions the energy is rapidly radiated away 
 during the collision of kinks with the anti\-kink (see the inset) and tends slowly to $\mathcal{E}=4$ for $s\ra \infty$. }}
	\label{fig:N3energy}
\end{figure}

\begin{figure}[ht]
	\centering
	\includegraphics[width=0.61\textwidth]{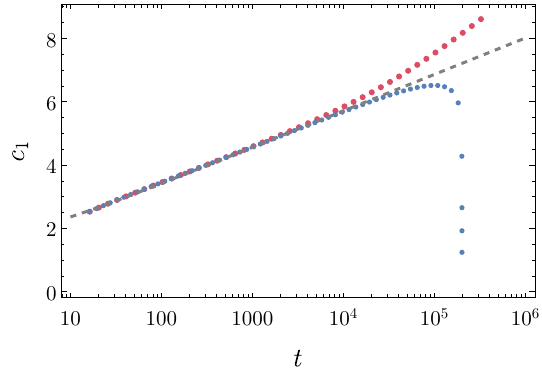}
	\caption{\small{Position of the kink $c_1(t)$ for sub- and supercritical solutions from Fig.~\ref{fig:N3evolution}. The gray dashed curve is the fit of formula \eqref{fitN3} to the numerical data for intermediate times. For later times the solutions deviate from criticality along the unstable directions. }
	}
	\label{fig:N3location}
\end{figure}

The snapshots of marginally sub- and supercritical solutions are depicted in Fig.~\ref{fig:N3evolution}, while 
the corresponding Bondi energies and positions of the kink are shown in Fig.~\ref{fig:N3energy} and Fig.~\ref{fig:N3location}.
\vskip 0.1cm
\noindent\emph{Acknowledgement.} 
 This work was  supported by the National Science Centre grant No.~2017/26/A/ST2/00530 (P.B.), the ERC Starting Grant \mbox{``INSOLIT''} 101117126 (J.J.), the Austrian Science Fund (FWF) Project \href{http://doi.org/10.55776/P36455}{P 36455} and START-Project \href{http://doi.org/10.55776/Y963}{Y 963} (M.M.). M.M. acknowledges the hospitality and support of Laboratoire Analyse, G\'eom\'etrie et Applications, Sorbonne Paris North University, during the research stay when part of this work was conducted.

\end{document}